\documentclass[12pt,a4paper]{article}

\usepackage[margin=1in]{geometry}
\usepackage{mathptmx}                 
\usepackage[T1]{fontenc}
\usepackage{amsmath,amssymb,amsthm}
\usepackage{graphicx}
\usepackage{booktabs,multirow}
\usepackage[font=small,labelfont=bf]{caption}
\usepackage{subcaption}
\usepackage[numbers,sort]{natbib}
\usepackage{xcolor}
\usepackage{float}
\usepackage{microtype}
\usepackage{authblk}
\usepackage{algorithm}
\usepackage{algpseudocode}
\usepackage{tikz}
\usetikzlibrary{positioning,arrows.meta,calc}
\usepackage[colorlinks=true,linkcolor=blue!50!black,citecolor=blue!50!black,urlcolor=blue!50!black]{hyperref}
\usepackage[capitalise,nameinlink]{cleveref}
\usepackage{authblk}

\graphicspath{{figures/}}

\newcommand{\bestLogistic}{Chebyshev KAN}

\newcommand{\bestSingularIVP}{B-spline KAN}

\newcommand{\bestSingularBVP}{B-spline KAN}

\newcommand{\bestFalknerSkan}{B-spline KAN}

\newcommand{\bestFourthOrder}{Chebyshev KAN}
\newcommand{\gridGainLogistic}{326.0}
\newcommand{\gridFinalLogistic}{48}
\newcommand{\gridGainSingularIVP}{3.8}

\newcommand{\gridGainSingularBVP}{0.7}

\newcommand{\gridGainFalknerSkan}{4.9}

\newcommand{\gridGainFourthOrder}{1.0}

\newcommand{\paramsKANmin}{80}
\newcommand{\paramsKANmax}{245}
\newcommand{\paramsMLP}{8{,}513}
\newcommand{\nUnfit}{8}
\newcommand{\nRuns}{225}

\providecommand{\gridGainLogistic}{??}\providecommand{\gridGainSingularIVP}{??}
\providecommand{\gridGainSingularBVP}{??}\providecommand{\gridGainFalknerSkan}{??}
\providecommand{\gridGainFourthOrder}{??}\providecommand{\gridFinalLogistic}{??}

\newif\ifshowchecks \showchecksfalse

\newtheorem{remark}{Remark}
\crefname{definition}{Definition}{Definitions}
\crefname{remark}{Remark}{Remarks}
\crefname{algorithm}{Algorithm}{Algorithms}

\newcommand{\params}[2]{\csname params#1#2\endcsname}  

\title{Neural Network Solvers for Nonlinear ODEs: Benchmarking MLPs and B-Spline/Chebyshev KANs}
\author[1]{Mahidhar Reddy Patukuri}
\author[1]{B. Veena S. N. Rao}
\author[1,$\dagger$]{S. M. Mallikarjunaiah}
\affil[1]{Department of Mathematics \& Statistics, Texas A\&M University-Corpus Christi, Texas- 78412, USA}
\affil[$\dagger$]{Corresponding author}
\affil[ ]{\textit{E-mail addresses:} \texttt{mpatukuri@islander.tamucc.edu} (Mahidhar Reddy Patukuri), \texttt{M.Muddamallappa@tamucc.edu} (S.M. Mallikarjunaiah)
\texttt{bv.rao@tamucc.edu} (B. Veena S. N. Rao)}
\date{}

\begin{document}
\maketitle

\begin{abstract}
\noindent
Three neural networks are compared as solvers of nonlinear ordinary differential equations: a multilayer perceptron (MLP), a Kolmogorov-Arnold network (KAN) whose edges are B-splines, and a KAN whose edges are Chebyshev series. Each network is trained only on the residual of the equation and on its initial or boundary conditions. The exact solution, or a numerical solution accurate to $10^{-13}$, is used only afterward to measure the error. Five test problems of differential orders one to four are investigated: the logistic equation, a singular Lane-Emden type initial value problem, a singular boundary value problem, the Falkner-Skan equation with suction, a fourth-order problem. The same training budget and the same choice of five activation functions are allocated to all networks. No numerical failures are encountered in any of the \nRuns{} training runs. For every problem, the best KAN is found to be at least as accurate as the best MLP while \paramsKANmin{} to \paramsKANmax{} parameters are used instead of \paramsMLP{}; the margin ranges from a tie to a factor of three. In the Falkner-Skan setting of an earlier study, the B-spline KAN is shown to be 34 to 220 times more accurate at every published test point. Much smaller MLPs are also found to be competitive, whereas at orders three and four, failures are experienced more frequently by KANs, and a coarse basis is required. The derivatives in the equation are distorted by the usual grid extension of a B-spline KAN; consequently, it is replaced with knot insertion, by which the logistic error is lowered more than a hundredfold, though no improvement is provided at order four.
\end{abstract}

\noindent\textbf{Keywords}\quad Physics-informed neural networks $\cdot$ Kolmogorov-Arnold networks
$\cdot$ Chebyshev polynomials $\cdot$ B-splines $\cdot$ Falkner-Skan equation $\cdot$ Singular
boundary value problems $\cdot$ Grid extension

\section{Introduction}\label{sec:intro}

Many nonlinear ordinary differential equations (ODEs) arising in physics and engineering are singular, posed on truncated semi-infinite domains, or of high order, and typically lack closed-form solutions. Classical numerical methods discretize the domain to solve a finite-dimensional surrogate; however, higher-order nonlinear problems often require complex techniques such as shooting, linearization, or specially designed grids. To overcome these limitations, neural-network solvers replace traditional domain discretization with smooth, differentiable parametric functions trained to satisfy the governing differential equation at collocation points alongside initial or boundary data \citep{lee1990, meade1994, lagaris1998, raissi2019}. Grounded in the universal approximation theorem for feedforward networks \citep{cybenko1989, hornik1989}, this approach has evolved into the broad field of physics-informed machine learning \citep{karniadakis2021, cuomo2022}, encompassing strong-form formulations \citep{raissi2019, sirignano2018}, variational approaches \citep{e2018ritz}, and specialized software frameworks \citep{lu2021deepxde}. For instance, feedforward networks trained via mean-squared residual loss have successfully reproduced complex models such as the Thomas-Fermi equation and Falkner-Skan nanofluid flows without explicit boundary condition adjustments \citep{venkatachalapathy2023, yacob2011, martinez2024approximation, venkatachalapathy2023deep, mallikarjunaiah2023deep, martinez2026numerical}.

Despite their versatility, training physics-informed neural networks (PINNs) remains notoriously challenging. Wang et al.\ \citep{wang2021gradient} identified gradient imbalance between residual and boundary terms as a primary failure mode, prompting adaptive loss weighting schemes, while neural-tangent-kernel analyses explain disparities in component convergence rates \citep{wang2022ntk}. Furthermore, PINNs can fail on basic convection-reaction problems \citep{krishnapriyan2021}, and standard feedforward architectures inherently learn low-frequency components first \citep{rahaman2019}. While remedies like adaptive activation functions \citep{jagtap2020} and point-wise self-adaptive loss weights \citep{mcclenny2023} have been proposed, the present study intentionally adopts a fixed, standardized protocol across all tested architectures.

Historically, nearly all physics-informed ODE solvers rely on multilayer perceptrons (MLPs), which feature fixed node nonlinearities and learnable scalar weights on every edge. Recently, Kolmogorov-Arnold Networks (KANs) \citep{liu2024kan}, inspired by the Kolmogorov-Arnold representation theorem \citep{kolmogorov1957, arnold1957, schmidthieber2021}, introduced a paradigm shift by moving learnable nonlinearities onto the edges themselves. Each edge carries an individual univariate function parameterized either by B-splines \citep{liu2024kan} or orthogonal polynomial expansions such as Chebyshev series \citep{ss2024chebyshev, chebpikan2025}. KANs claim several advantages for scientific computing: high parameter efficiency, interpretable edge functions that can be mapped to closed-form expressions \citep{liu2025kan2}, and the ability to refine B-spline KANs post-training via grid extension.

However, it remains unclear whether these purported benefits persist when networks are trained directly on differential equation residuals. Because the residual of a $p$-th order ODE differentiates the network $p$ times, and a degree-$k$ B-spline possesses only $k-1$ continuous derivatives, numerical robustness can be severely tested. Moreover, existing literature presents conflicting conclusions: Shukla et al.\ \citep{shukla2024comprehensive} found B-spline KANs less efficient than MLPs (though polynomial variants were competitive), whereas Wang et al.\ \citep{wang2024kinn} reported superior accuracy for KAN-based solvers in solid mechanics. Crucially, none of these existing studies systematically evaluate performance across varying differential orders.

To bridge this gap, this paper provides a comprehensive evaluation with the following main contributions:
\begin{enumerate}
    \item We compare an MLP, a B-spline KAN, and a Chebyshev KAN across five nonlinear ODEs of orders one to four—including two singular problems and a new fourth-order problem with an exact closed-form solution—under an identical training budget (\cref{sec:methods,sec:problems}).
    \item We employ a robust training procedure where all \nRuns{} runs achieve numerical stability, documenting every instance where a run fails to capture the exact solution (\cref{tab:robust}).
    \item We demonstrate that KAN bases must begin with a coarse resolution at higher orders three and four (\cref{rem:scale}), and we replace standard grid extension with derivative-preserving knot insertion (\cref{rem:nested}, \cref{sec:gridext}).
    \item We extract and plot every learned edge function from the trained KANs, matching them directly to simple analytic formulas (\cref{sec:interp}).
\end{enumerate}

\section{Methodology}\label{sec:methods}

\subsection{Physics-informed loss}\label{sec:loss}
Every problem in this paper has the form
\begin{equation}
\mathcal N[u](x) = 0 \ \text{ for } x\in(a,b), \qquad \mathcal B_j[u] = g_j \ \text{ for } j = 1,\dots,J.
\label{eq:ode}
\end{equation}
Here $u$ is the unknown solution on the interval $[a,b]$, and $\mathcal N[u]$ is the \emph{residual} of the
equation: the left-hand side minus the right-hand side, which is zero exactly when $u$ solves the equation. For the
logistic equation $y'=y(1-y)$, for example, $\mathcal N[u](x)=u'(x)-u(x)\,(1-u(x))$. The order $p$ of the equation is
the highest derivative in $\mathcal N$. The $J$ conditions fix the initial or boundary values: each $\mathcal B_j$
picks out the value of $u$ or of one of its derivatives at an end point, and $g_j$ is the prescribed value. The
logistic problem has one condition, $u(0)=\tfrac12$, so $J=1$, $\mathcal B_1[u]=u(0)$ and $g_1=\tfrac12$. The five
problems have $J=1,2,2,3,4$ conditions (\cref{tab:problems}).

The network $u_\theta$, with trainable parameters $\theta$, replaces $u$. Training minimises
\begin{equation}
\begin{aligned}
\mathcal L(\theta) &= \mathcal L_{\mathrm{res}}(\theta) + \lambda\,\mathcal L_{\mathrm{bc}}(\theta),\\
\mathcal L_{\mathrm{res}}(\theta) &= \frac{1}{N_r}\sum_{i=1}^{N_r}\big(\mathcal N[u_\theta](x_i)\big)^2,\qquad
\mathcal L_{\mathrm{bc}}(\theta) = \sum_{j=1}^{J}\big(\mathcal B_j[u_\theta]-g_j\big)^2.
\end{aligned}
\label{eq:loss}
\end{equation}
The terms are as follows.
\begin{itemize}
\item $x_1,\dots,x_{N_r}$ are the \emph{collocation points}: $N_r=1{,}000$ points drawn uniformly at random in
$(a,b)$. During the Adam phase a new set is drawn at every step; the L-BFGS phase uses one fixed set
(\cref{sec:training}).
\item $\mathcal L_{\mathrm{res}}$ is the mean squared residual at the collocation points. It is zero only if the
network satisfies the differential equation at every one of them.
\item $\mathcal L_{\mathrm{bc}}$ is the sum of the squared errors in the $J$ conditions. It is zero only if the
network has the prescribed initial or boundary values.
\item $\lambda>0$ is the weight of the condition term relative to the residual term. The residual alone is zero for
every solution of the equation, so it is the conditions that single out the correct one. $\mathcal L_{\mathrm{res}}$
averages over 1{,}000 points while $\mathcal L_{\mathrm{bc}}$ has only $J\le4$ terms, and with a small $\lambda$ an
optimizer can reduce the residual while the conditions stay slightly wrong, which shifts the whole solution.
\end{itemize}
We use $\lambda=100$ for every problem and every architecture. The value was chosen in the pilot study of
\cref{sec:pilot} by comparing $\lambda=10$ and $\lambda=100$ on the training loss of the fourth-order problem, the
one with four conditions at a single end point. With $\lambda=100$ all three networks reached a lower training loss,
and their errors dropped accordingly (MLP from $1.7\times10^{-3}$ to $4.1\times10^{-6}$, Chebyshev KAN from
$3.5\times10^{-2}$ to $3.6\times10^{-4}$, B-spline KAN from $0.45$ to $0.09$ in that pilot run, which was much shorter
than the final training). We did not tune
$\lambda$ further or per problem, so that no architecture benefits from its own choice; adaptive weights
\citep{wang2021gradient,mcclenny2023} are an alternative we did not use.

All derivatives in $\mathcal N[u_\theta]$ and $\mathcal B_j[u_\theta]$ are computed by automatic differentiation
\citep{baydin2018} in PyTorch \citep{paszke2019}, exactly and without finite differences. No values of the solution
enter \cref{eq:loss}.

\subsection{Input normalisation}
Each network receives $z = 2(x-a)/(b-a)-1\in[-1,1]$ rather than $x$. This matters for KANs, whose bases
are defined on a fixed interval: without it, most of the domain $[0,3]$ of the Falkner-Skan problem would
lie outside the support of a spline grid on $[-1,1]$. The map is applied identically to all three
architectures, and derivatives are still taken with respect to $x$.

\subsection{Multilayer perceptron}
\begin{figure}[H]
\centering
\includegraphics[width=0.45\linewidth]{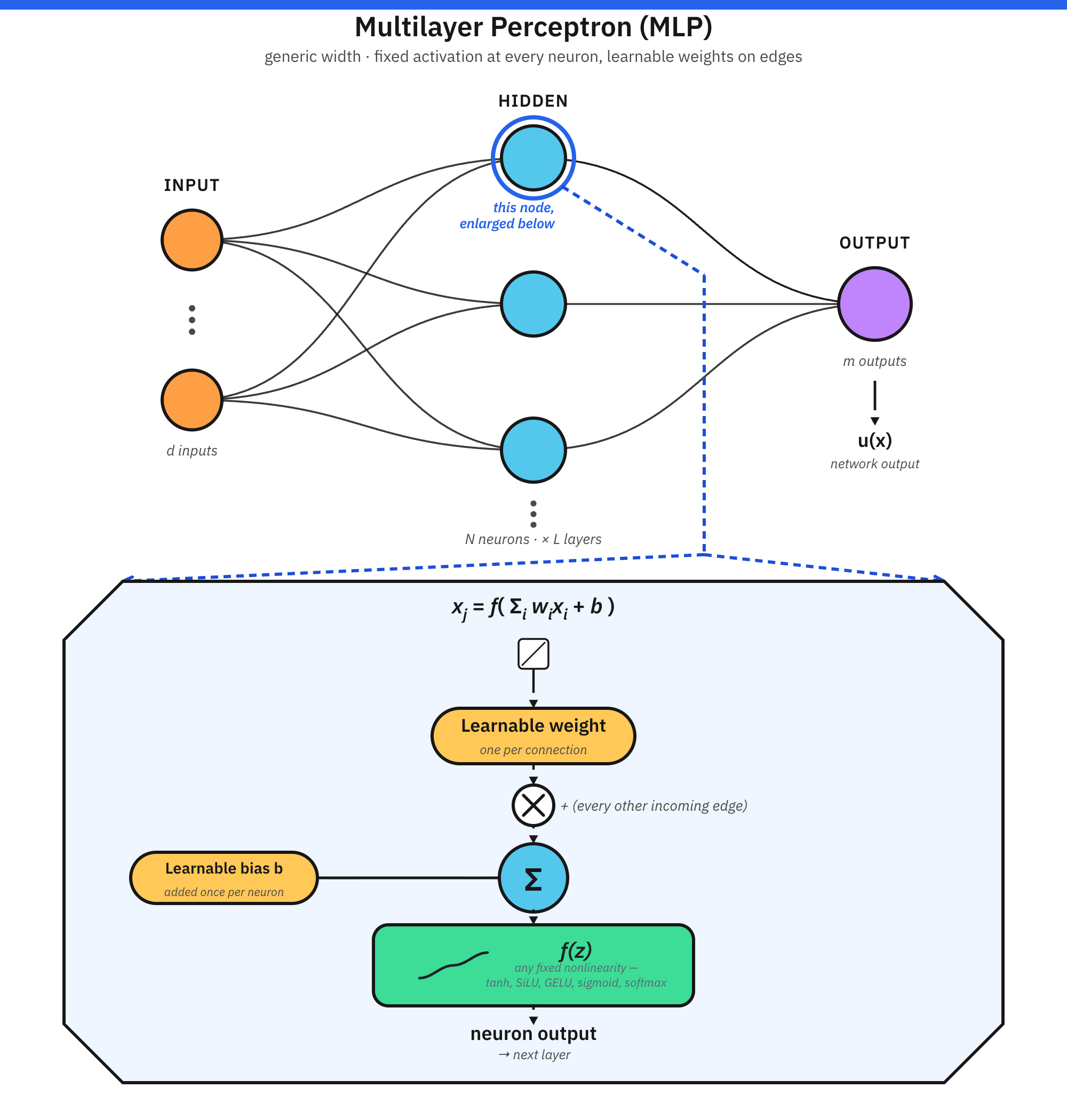}
\caption{Multilayer perceptron. Every edge carries a learnable scalar weight; every neuron adds a learnable
bias to the weighted sum of its inputs and applies the same fixed activation $f$. In this study
$d=m=1$ and there are three hidden layers of 64 neurons, $f\in\{\tanh,\mathrm{SiLU},\mathrm{GELU},
\mathrm{sigmoid},\mathrm{softmax}\}$.}
\label{fig:arch-mlp}
\end{figure}

The baseline is the standard PINN MLP with three hidden layers of 64 units,
\begin{equation}
u_\theta(z) = W_4\,\sigma\!\big(W_3\,\sigma(W_2\,\sigma(W_1 z+b_1)+b_2)+b_3\big)+b_4,
\label{eq:mlp}
\end{equation}
with Glorot-normal initial weights \citep{glorot2010} and \params{Logistic}{MLP} trainable parameters.
The activation $\sigma$ is one of $\tanh$, SiLU \citep{elfwing2018}, GELU \citep{hendrycks2016}, sigmoid and
softmax (applied across the units of a layer), the last two being among those studied in \citep{venkatachalapathy2023}. \Cref{fig:arch-mlp} shows
the structure of a single neuron.

\subsection{B-spline Kolmogorov-Arnold network}
\begin{figure}[H]
\centering
\includegraphics[width=0.45\linewidth]{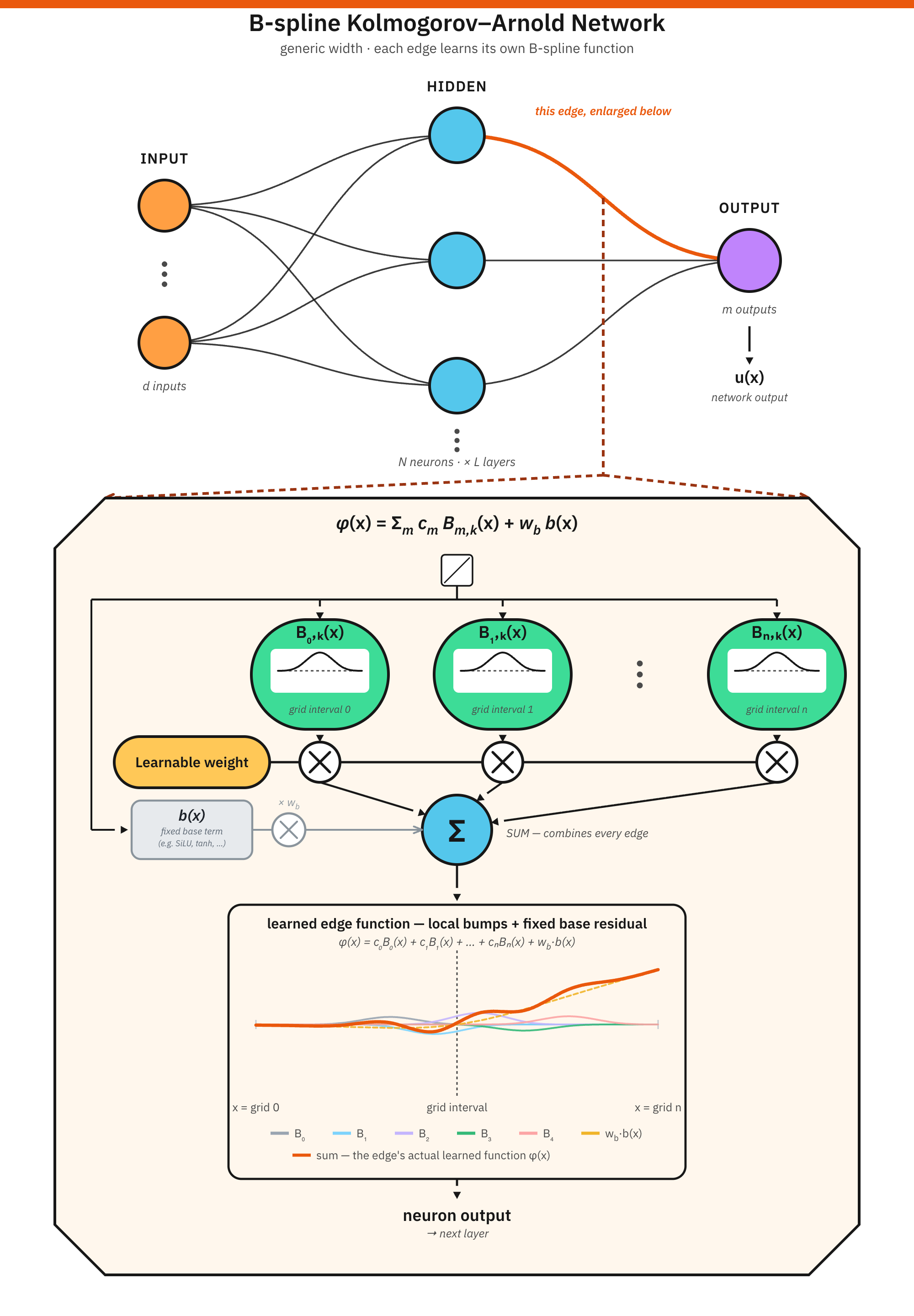}
\caption{B-spline Kolmogorov-Arnold network. Every edge learns its own function $\phi_{ij}$ as a weighted
sum of B-spline basis functions on a grid plus a residual base term (\cref{eq:bspline}); nodes only add
their incoming edges. In this study $d=m=1$, the width is $[1,5,1]$, the grid has $G=3$ intervals, the
spline degree is $k=\max(3,p+1)$, and the base function is selected from SiLU (shown), $\tanh$, GELU,
sigmoid or none.}
\label{fig:arch-bspline}
\end{figure}

A KAN layer with $n_{\mathrm{in}}$ inputs and $n_{\mathrm{out}}$ outputs computes
$y_j=\sum_{i=1}^{n_{\mathrm{in}}}\phi_{ij}(x_i)$, where every edge carries its own learnable function.
Following \citep{liu2024kan},
\begin{equation}
\phi_{ij}(x) = w^{b}_{ij}\,b(x) + w^{s}_{ij}\sum_{m=1}^{G+k} c_{ij,m}\,B_{m,k}(x),
\label{eq:bspline}
\end{equation}
where $B_{m,k}$ are B-splines of degree $k$ \citep{deboor2001} on $G$ uniform intervals (extended by $k$ knots on each side)
and $b$ is a residual base function; Liu et al.\ \citep{liu2024kan} use $b=\mathrm{SiLU}$. Mirroring the
MLP, $b$ is one of SiLU, $\tanh$, GELU and sigmoid, or the base branch is omitted (``none'').
We use the width $[1,5,1]$ with $G=3$, i.e.\ \params{FalknerSkan}{BSpline} parameters for the
Falkner-Skan problem. \Cref{fig:arch-bspline} shows how a single edge function is built.

\begin{remark}[Spline degree and differential order]\label{rem:degree}
The $k$-th derivative of a degree-$k$ spline is piecewise constant and its $(k{+}1)$-th derivative
vanishes almost everywhere. A residual of order $p$ therefore needs $k\ge p+1$ for the $p$-th
derivative to be continuous. We set $k=\max(3,p+1)$, i.e.\ $k=3,3,3,4,5$ for P1 to P5.
\end{remark}

\begin{remark}[Grid resolution and derivative scale]\label{rem:scale}
On a grid of spacing $h$ the $p$-th derivative of a B-spline scales like $h^{-p}$. For $p=3,4$ a
fine grid therefore produces residual gradients that are several orders of magnitude larger than those of
a coarse grid, and the optimization stalls. In a pilot study (training loss only, \cref{sec:pilot})
$G=10$ left the Falkner-Skan loss on the trivial plateau, while $G=5$ and $G=3$ converged, and on the fourth-order problem $G=3$ was clearly better than $G=5$. We therefore train on a
coarse grid and add resolution by grid extension (\cref{sec:gridext}).
\end{remark}

As in \citep{liu2024kan}, the grids of hidden layers are adapted to the distribution of their inputs during
the first 30\% of the Adam iterations; the spline coefficients are re-fitted by least squares at every
update, so that the values of the network are preserved.

\begin{remark}[Derivative-exact grid extension]\label{rem:nested}
Grid extension \citep{liu2024kan} replaces the grid of a trained B-spline KAN by a finer one and fits the new
coefficients to the old spline by least squares. If the new knots are placed independently of the old ones,
the fit preserves the values of $\phi_{ij}$ only approximately and its $p$-th derivative much less so, and
a $p$-th order residual jumps after the extension. If instead every old knot is kept
and each interval is subdivided (knot insertion \citep{boehm1980,cohen1980}), the old spline space is a subspace of the new one, so the
old spline is represented \emph{exactly}, with all its derivatives, and the physics loss is unchanged by
the extension. We use this nested refinement, doubling the number of intervals on the knot span covered by
the data at each extension.
\end{remark}

\subsection{Chebyshev Kolmogorov-Arnold network}
\begin{figure}[H]
\centering
\includegraphics[width=0.45\linewidth]{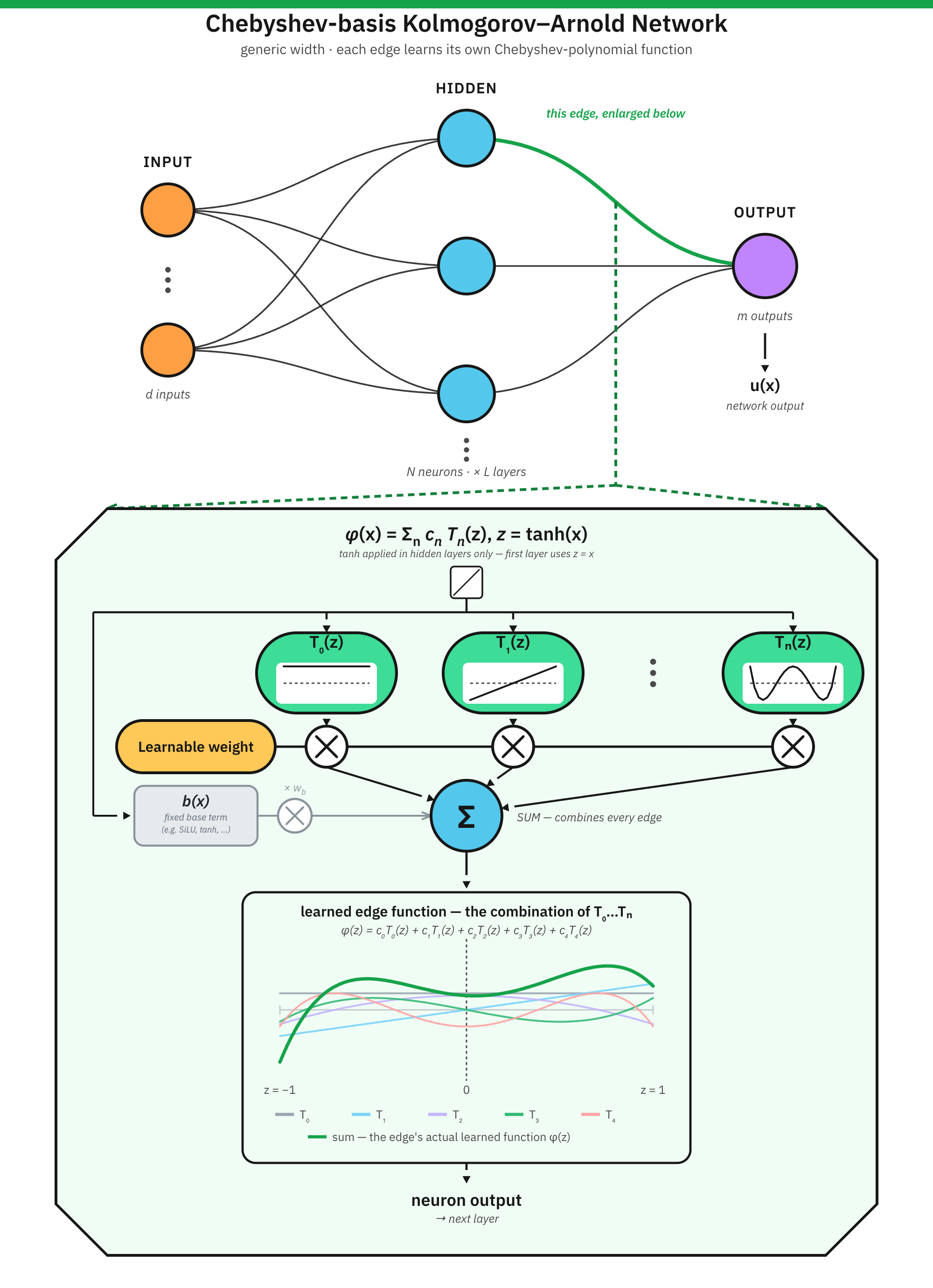}
\caption{Chebyshev Kolmogorov-Arnold network. Every edge learns its own function $\phi_{ij}$ as a weighted
sum of Chebyshev polynomials $T_0,\dots,T_d$ (\cref{eq:cheby}); nodes only add their incoming edges. In this
study $d_{\mathrm{in}}=m=1$, the width is $[1,5,5,1]$ and the degree is $d=5$; the first layer evaluates
$T_n$ directly on the normalised input, hidden layers on $\tanh$ of their input, and an optional residual
base term $w^{b}_{ij}\,b(x)$ is added (\cref{eq:cheby}).}
\label{fig:arch-cheby}
\end{figure}

The Chebyshev KAN replaces \cref{eq:bspline} by a truncated Chebyshev series
\begin{equation}
\phi_{ij}(x) = w^{b}_{ij}\,b(x) + \sum_{n=0}^{d} c_{ij,n}\,T_n\big(s(x)\big),\qquad T_{n+1}=2sT_n-T_{n-1},
\label{eq:cheby}
\end{equation}
with $s(x)=x$ in the first layer (whose input is already in $[-1,1]$) and $s(x)=\tanh x$ in the hidden
layers. Evaluating $T_n$ by the three-term recurrence \citep{mason2003,trefethen2019} keeps every derivative bounded on $[-1,1]$; the
alternative $\cos(n\arccos s)$ has an unbounded derivative at $s=\pm1$. The initial coefficients are drawn
with a standard deviation proportional to $(n+1)^{-2}$ so that the network starts smooth. The residual base
function $b$ is chosen from the same five options as for the B-spline KAN. We use the width
$[1,5,5,1]$ and degree $d=5$ (\params{FalknerSkan}{Cheby} parameters); \cref{fig:arch-cheby} shows a single
edge.

\subsection{Training protocol}\label{sec:training}
\Cref{alg:training} is used for every architecture and problem. It combines Adam \citep{kingma2015} with a
subsequent L-BFGS phase \citep{liu1989lbfgs}, the usual two-stage optimisation for PINNs \citep{lu2021deepxde}. The budget is identical for all three
networks (\cref{tab:hyper}); only the network inside $u_\theta$ changes.

\begin{algorithm}[H]
\caption{Physics-informed training (identical for all architectures)}\label{alg:training}
\begin{algorithmic}[1]
\Require network $u_\theta$, problem (\cref{eq:ode}), seed $s$, budget of \cref{tab:hyper}
\State set seed $s$ (initial weights and collocation sampling); draw a fixed set $X_f$ of $N_f$ interior points
\For{$t=1,\dots,T_{\mathrm{Adam}}$}
  \State sample $N_r$ interior points uniformly; evaluate $\mathcal L$ (\cref{eq:loss})
  \If{$\mathcal L$ or $\nabla_\theta\mathcal L$ is not finite} \State roll back to last checkpoint; halve the learning rate; \textbf{continue} \EndIf
  \State clip $\|\nabla_\theta\mathcal L\|\le 1$; Adam step; decay learning rate exponentially
  \State every 100 steps: evaluate $\mathcal L$ on $X_f$; keep the parameters with the lowest value
\EndFor
\For{$t=1,\dots,T_{\mathrm{LBFGS}}$} \Comment{on the fixed set $X_f$}
  \State L-BFGS step (strong Wolfe line search); stop and roll back if not finite; update best parameters
\EndFor
\State \Return parameters with the lowest $\mathcal L$ on $X_f$ \Comment{no solution data used}
\end{algorithmic}
\end{algorithm}

\begin{table}[H]
\centering
\caption{Training budget shared by all architectures and problems.}
\label{tab:hyper}
\small
\begin{tabular}{ll}
\toprule
Precision & float64 \\
Interior points per Adam step $N_r$ / fixed set $N_f$ & 1{,}000 / 1{,}000 \\
Adam iterations, learning rate & 10{,}000,\; $10^{-2}\to10^{-5}$ (exponential) \\
Gradient clipping & global norm 1 \\
L-BFGS iterations & 5{,}000 (history 50, strong Wolfe) \\
Condition weight $\lambda$ & 100 \\
Activation choices & 5 per architecture \\
Random seeds & 3 per activation (0, 1, 2) \\
Test points & 2{,}001, uniform \\
\bottomrule
\end{tabular}
\end{table}

\subsection{Testing procedure and seed selection}\label{sec:metrics}
Each trained network is tested in the same way.
\begin{enumerate}
\item The network is evaluated at $M=2{,}001$ equally spaced test points $x_m=a+(m-1)(b-a)/2000$,
$m=1,\dots,M$, which include both end points. These are different from the random collocation points used in
training.
\item The network values $\hat y_m=u_\theta(x_m)$ are compared with the true solution $y_m=y(x_m)$. For P1, P2 and P5
$y$ is the closed-form solution of \cref{tab:problems}. P3 and P4 have no closed form, so $y$ is a numerical
solution computed with \texttt{solve\_bvp} from SciPy \citep{virtanen2020}. This solver approximates the solution by a
piecewise cubic function and refines its mesh until the residual of the equation, relative to the size of the
solution, is below a tolerance \citep{kierzenka2001}; we use $10^{-10}$. To check the accuracy of these reference
solutions, we solved both problems again with a hundred times smaller tolerance ($10^{-12}$) on a different initial
mesh. At the test points, the two solutions differ by at most $2\times10^{-14}$ for P3 and $2\times10^{-13}$ for P4.
For P4, the reference also agrees to $10^{-8}$ with the Runge-Kutta values published in
\citep{venkatachalapathy2023}, which is all the digits printed there. Since every network error reported in this
paper is $10^{-7}$ or larger, the error of the reference solutions has no visible effect, and we treat them as exact.
\item From the $M$ differences we compute three error measures:
\[
\mathrm{rel}L^2=\frac{\big(\sum_{m}(\hat y_m-y_m)^2\big)^{1/2}}{\big(\sum_{m}y_m^2\big)^{1/2}},\qquad
\mathrm{MAE}=\frac1M\sum_{m}\lvert\hat y_m-y_m\rvert,\qquad
\max_m\lvert\hat y_m-y_m\rvert .
\]
The relative $L^2$ error is the main measure; it is independent of the scale of the solution, so it can be compared
across problems.
\end{enumerate} Each architecture is trained with five activation choices and three seeds per activation;
the main results use, for every problem and architecture, the run (activation and seed) with the lowest final
training loss on $X_f$, a choice that requires no knowledge of the solution. The activation ablation is reported
in \cref{sec:activations} and the behaviour of all runs in \cref{tab:robust}.

\subsection{Pilot study}\label{sec:pilot}
The learning rate and the KAN sizes were fixed by a short pilot run (1{,}500 Adam iterations,
seed 0) on the two highest-order problems, judged by training loss only. A learning rate of $10^{-2}$
was better than $10^{-3}$ for all three architectures. The condition weight was compared in a second, longer
pilot run on P5 (3{,}000 Adam and 1{,}000 L-BFGS iterations, seed 0), where $\lambda=100$ gave a lower training loss
than $\lambda=10$ for all three architectures (\cref{sec:loss}). For both KANs, coarse bases (B-spline $G\le5$, best $G=3$,
Chebyshev $d\le5$) trained well on P4 and P5 whereas $G=10$ and $d=8$ stalled (\cref{rem:scale}).

\section{Benchmark problems}\label{sec:problems}
\Cref{tab:problems} lists the five problems in order of differential order. For the singular problems
the residual is multiplied by $x$, which removes the division from the computational graph without
changing the equation on $(0,1]$.

\begin{table}[H]
\centering
\caption{Benchmark problems. The loss of each is \cref{eq:loss} with the residual and conditions listed.
$^{*}$No closed form. The test solution is computed with \texttt{solve\_bvp} and is accurate to better than
$10^{-12}$ at all 2{,}001 test points (\cref{sec:metrics}). Check values: P3, $y(0)=0.954135$; P4,
$f(3)=3.308343$, which agrees with the value $3.308343141$ published in \citep{venkatachalapathy2023}.}
\label{tab:problems}
\footnotesize
\setlength{\tabcolsep}{3pt}
\begin{tabular}{@{}llll@{}}
\toprule
 & Equation and conditions & Domain & Test solution \\
\midrule
P1 & $y'=y(1-y)$,\; $y(0)=\tfrac12$ & $[0,5]$ & exact $1/(1+e^{-x})$ \\
P2 & $xy''+y'-xy^3+3xy^5=0$,\; $y(0)=1,\ y'(0)=0$ & $[0,1]$ & exact $(1+x^2)^{-1/2}$ \\
P3 & $xy''+3y'=\tfrac{x}{2}-\tfrac{x}{8y^2}$,\; $y'(0)=0,\ y(1)=1$ & $[0,1]$ & no closed form; numerical$^{*}$ \\
P4 & $f'''+ff''+\beta(1-f'^2)=0$,\; $f(0)=f_w,\ f'(0)=0,\ f'(3)=1$ & $[0,3]$ & no closed form; numerical$^{*}$ \\
P5 & $y^{(4)}=24y^5$,\; $y(0)=1,\ y'(0)=-1,\ y''(0)=2,\ y'''(0)=-6$ & $[0,1]$ & exact $1/(1+x)$ \\
\bottomrule
\end{tabular}
\end{table}

\paragraph{P1: logistic growth ($p=1$).} Verhulst's population model \citep{verhulst1838}; the simplest case:
first order, no singularity, closed form.

\paragraph{P2: singular IVP ($p=2$).} A Lane-Emden type equation \citep{wazwaz2001} $y''+\frac1x y'-y^3+3y^5=0$ with a
coordinate singularity at $x=0$ and the closed-form solution $y=(1+x^2)^{-1/2}$.

\paragraph{P3: singular BVP ($p=2$).} $y''+\frac3x y'=\frac12-\frac{1}{8y^2}$ with $y'(0)=0$, $y(1)=1$
has a rational nonlinearity and no closed form. Its test solution is computed with \texttt{solve\_bvp}
as described in \cref{sec:metrics}, using the solver's option for the singular term $3y'/x$. Inside the residual only, $y$ is floored at $0.1$
in the term $1/y^2$ so that an early excursion of the network towards $y=0$ cannot produce an infinite
loss; the solution satisfies $0.95\le y\le1$ and never reaches the floor.

\paragraph{P4: Falkner-Skan flow with suction ($p=3$).} The Falkner-Skan equation
\citep{falknerskan1931,hartree1937} describes the similarity boundary layer on a wedge; the suction and wedge
parameters are those of the nanofluid study \citep{yacob2011,venkatachalapathy2023}. With wedge parameter $m=0.0909$,
$\beta=2m/(m+1)$, suction $s=0.5$ and $f_w=2s/(m+1)$, on the truncated domain $\eta\in[0,3]$ used in
\citep{venkatachalapathy2023}. Its test solution is also computed with \texttt{solve\_bvp}
(\cref{sec:metrics}).

\paragraph{P5: nonlinear fourth-order IVP ($p=4$).} We introduce $y^{(4)}=24y^5$ with
$y(0)=1$, $y'(0)=-1$, $y''(0)=2$, $y'''(0)=-6$. Since $\frac{d^4}{dx^4}(1+x)^{-1}=24(1+x)^{-5}$, the exact
solution is $y=1/(1+x)$. It is the only problem whose residual contains a fourth derivative and whose
four conditions are all imposed at one end point, which makes it the hardest test of both the basis
smoothness (\cref{rem:degree}) and the derivative scale (\cref{rem:scale}).

\section{Results}\label{sec:results}

\subsection{Accuracy under an identical budget}\label{sec:accuracy}
\begin{table*}[t]
\centering
\caption{Best run (over activations and seeds, selected by training loss) of the three architectures under an identical training budget. Relative $L^2$ error and mean absolute error (MAE) are measured on 2{,}001 uniform test points against the exact (P1, P2, P5) or reference (P3, P4) solution. The best relative $L^2$ error per problem is in bold. Act.: selected activation (MLP: node activation; KANs: residual base function). $|\theta|$: trainable parameters.}
\label{tab:main}
\small
\setlength{\tabcolsep}{3pt}
\resizebox{\textwidth}{!}{%
\begin{tabular}{lcccrcccrcccr}
\toprule
 & \multicolumn{4}{c}{MLP} & \multicolumn{4}{c}{B-spline KAN} & \multicolumn{4}{c}{Chebyshev KAN} \\
\cmidrule(lr){2-5}\cmidrule(lr){6-9}\cmidrule(lr){10-13}
Problem & rel.\ $L^2$ & MAE & act. & $|\theta|$ & rel.\ $L^2$ & MAE & act. & $|\theta|$ & rel.\ $L^2$ & MAE & act. & $|\theta|$ \\
\midrule
Logistic ($p{=}1$) & $4.9\times10^{-5}$ & $3.7\times10^{-5}$ & tanh & 8{,}513 & $2.4\times10^{-5}$ & $1.6\times10^{-5}$ & sigm. & 80 & \textbf{$2.4\times10^{-5}$} & $1.8\times10^{-5}$ & none & 210 \\
Singular IVP ($p{=}2$) & $1.1\times10^{-5}$ & $8.7\times10^{-6}$ & GELU & 8{,}513 & \textbf{$1.1\times10^{-5}$} & $8.4\times10^{-6}$ & GELU & 80 & $2.5\times10^{-5}$ & $1.9\times10^{-5}$ & SiLU & 245 \\
Singular BVP ($p{=}2$) & $3.7\times10^{-7}$ & $2.7\times10^{-7}$ & SiLU & 8{,}513 & \textbf{$1.2\times10^{-7}$} & $8.3\times10^{-8}$ & sigm. & 80 & $4.7\times10^{-7}$ & $3.4\times10^{-7}$ & tanh & 245 \\
Falkner-Skan ($p{=}3$) & $1.6\times10^{-6}$ & $2.6\times10^{-6}$ & sigm. & 8{,}513 & \textbf{$5.4\times10^{-7}$} & $9.0\times10^{-7}$ & none & 80 & $4.8\times10^{-6}$ & $8.1\times10^{-6}$ & sigm. & 245 \\
Fourth-order IVP ($p{=}4$) & $1.6\times10^{-6}$ & $9.5\times10^{-7}$ & GELU & 8{,}513 & $2.6\times10^{-6}$ & $1.6\times10^{-6}$ & SiLU & 100 & \textbf{$7.0\times10^{-7}$} & $4.0\times10^{-7}$ & GELU & 245 \\
\bottomrule
\end{tabular}}
\end{table*}

\begin{figure}[H]
\centering
\includegraphics[width=\linewidth]{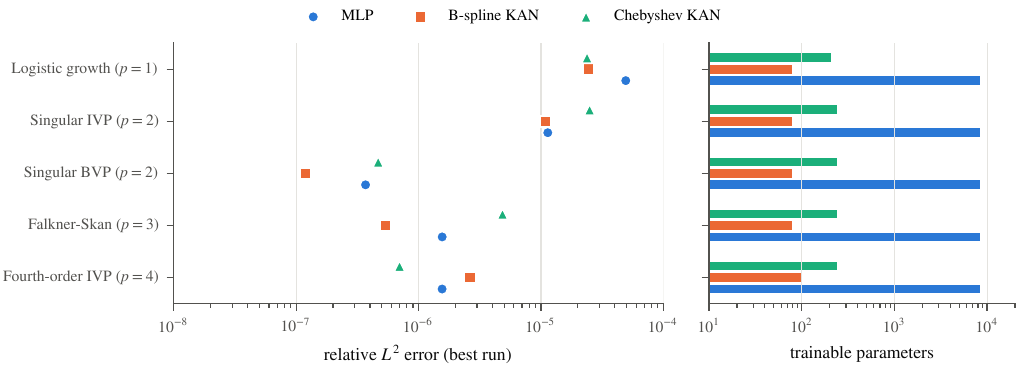}
\caption{Left: best-seed relative $L^2$ error per problem. Right: number of trainable parameters
(logarithmic). Markers and colors identify the architecture.}
\label{fig:summary}
\end{figure}

\begin{figure}[H]
\centering
\includegraphics[width=\linewidth]{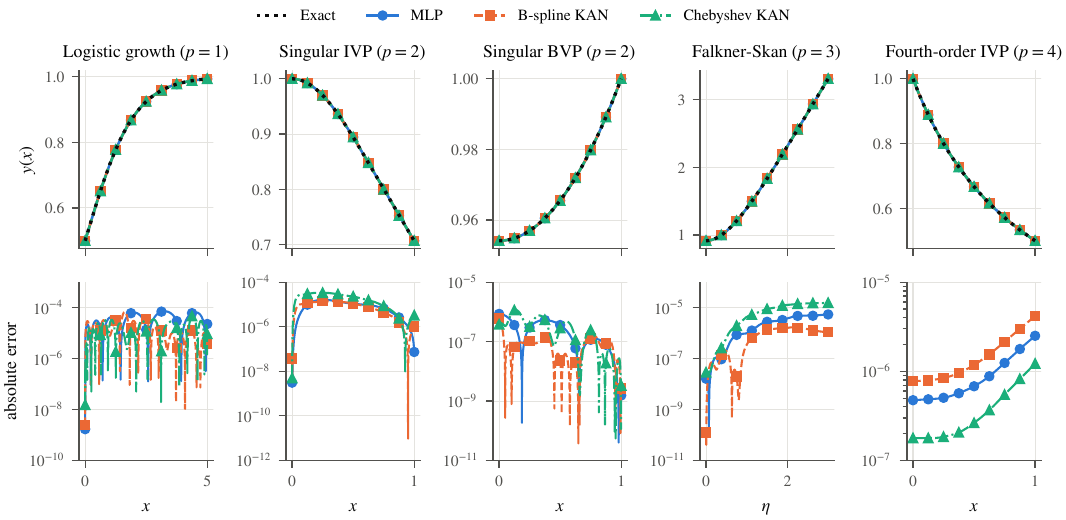}
\caption{Top: network solutions (best seed) against the exact or reference solution (dotted).
Bottom: pointwise absolute error, logarithmic scale.}
\label{fig:solutions}
\end{figure}

\Cref{tab:main} and \cref{fig:summary} summarise the best runs. The two KANs use between
\paramsKANmin{} and \paramsKANmax{} trainable parameters, one to two orders of magnitude fewer than the
\paramsMLP{} of the MLP. The lowest relative $L^2$ error is obtained by the \bestLogistic{} on P1, the
\bestSingularIVP{} on P2, the \bestSingularBVP{} on P3, the \bestFalknerSkan{} on P4 and the
\bestFourthOrder{} on P5. The margins differ: on P1 the best KAN is about twice as accurate as the MLP and on
P2 the two are practically tied, whereas on P3, P4 and P5 the best KAN is two to three times more accurate than the
MLP. \Cref{fig:solutions} shows that all three networks are visually indistinguishable from the exact or
reference solution. The pointwise errors differ in shape: for the initial value problems P1 and P5 the error
grows away from the initial point, for the Falkner-Skan problem it grows with $\eta$ towards the far-field
condition, and for the singular BVP P3 it is largest near the singular end point $x=0$. Pointwise errors at
equally spaced points of every domain are listed in \cref{tab:pointwise} (\cref{app:pointwise}).

\paragraph{Comparison with a published network solution.} The Falkner-Skan configuration is that of
\citep{venkatachalapathy2023}, who tabulate the pointwise error of their feedforward network (sigmoid
activation, AdaMax optimizer). \Cref{tab:fsref} places the errors of our best runs at the same points. At
every tabulated point, the error of the B-spline KAN (\params{FalknerSkan}{BSpline} parameters) is 34 to 220 times
smaller than the published one, that of the MLP 6 to 24 times smaller and that of the Chebyshev KAN 2 to 10 times
smaller; our reference solution agrees with theirs to $10^{-8}$, so the comparison is not limited by the
reference. Part of this gain comes from the training protocol, which all three networks share, and not from the
architecture alone.

\begin{table}[t]
\centering
\caption{Falkner-Skan problem ($m=0.0909$, $s=0.5$): absolute error $|f_{\mathrm{ref}}-\hat f|$ at the points tabulated by Venkatachalapathy and Mallikarjunaiah \citep{venkatachalapathy2023} (their Table 5, feedforward network with sigmoid activation), next to the best runs of this study.}
\label{tab:fsref}
\small
\setlength{\tabcolsep}{4pt}
\ifdefined\fitbox\else\newsavebox\fitbox\fi
\sbox\fitbox{\begin{tabular}{cccccc}
\toprule
$\eta$ & $f_{\mathrm{ref}}(\eta)$ & Ref.~\citep{venkatachalapathy2023} & MLP & B-spline KAN & Chebyshev KAN \\
\midrule
0.3 & 0.971145 & $2.2\times10^{-6}$ & $1.1\times10^{-7}$ & $5.5\times10^{-8}$ & $2.1\times10^{-7}$ \\
0.6 & 1.113635 & $1.1\times10^{-5}$ & $4.8\times10^{-7}$ & $5.0\times10^{-8}$ & $1.2\times10^{-6}$ \\
0.9 & 1.317738 & $6.0\times10^{-6}$ & $9.5\times10^{-7}$ & $1.2\times10^{-7}$ & $3.3\times10^{-6}$ \\
1.2 & 1.562807 & $2.8\times10^{-5}$ & $1.6\times10^{-6}$ & $7.6\times10^{-7}$ & $5.8\times10^{-6}$ \\
1.5 & 1.833419 & $5.6\times10^{-5}$ & $2.8\times10^{-6}$ & $1.4\times10^{-6}$ & $8.9\times10^{-6}$ \\
1.8 & 2.118838 & $5.9\times10^{-5}$ & $3.1\times10^{-6}$ & $1.6\times10^{-6}$ & $1.1\times10^{-5}$ \\
2.1 & 2.412195 & $5.8\times10^{-5}$ & $4.3\times10^{-6}$ & $1.7\times10^{-6}$ & $1.3\times10^{-5}$ \\
2.4 & 2.709473 & $8.4\times10^{-5}$ & $4.9\times10^{-6}$ & $1.5\times10^{-6}$ & $1.5\times10^{-5}$ \\
2.7 & 3.008534 & $1.2\times10^{-4}$ & $4.9\times10^{-6}$ & $1.2\times10^{-6}$ & $1.5\times10^{-5}$ \\
\bottomrule
\end{tabular}}
\ifdim\wd\fitbox>\linewidth\resizebox{\linewidth}{!}{\usebox\fitbox}\else\usebox\fitbox\fi
\end{table}

\subsection{Choice of activation}\label{sec:activations}
\begin{figure}[H]
\centering
\includegraphics[width=\linewidth]{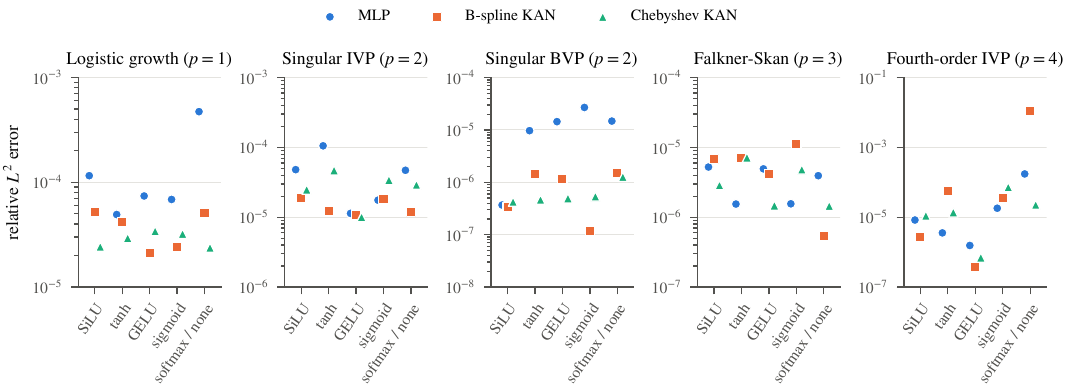}
\caption{Activation ablation: relative $L^2$ error of the best seed for each activation (MLP: node
activation; KANs: residual base function). The last category is softmax for the MLP and no base branch for
the KANs.}
\label{fig:activations}
\end{figure}

\begin{table}[t]
\centering
\caption{Activation ablation: relative $L^2$ error of the best seed for every activation. MLP: activation of the hidden nodes; KANs: residual base function $b$ of every edge. The last column is softmax for the MLP and no base branch (pure basis) for the KANs. Bold: the activation selected for \cref{tab:main}.}
\label{tab:activations}
\scriptsize
\setlength{\tabcolsep}{2.5pt}
\ifdefined\fitbox\else\newsavebox\fitbox\fi
\sbox\fitbox{\begin{tabular}{llccccc}
\toprule
Problem & Architecture & SiLU & tanh & GELU & sigmoid & softmax / none \\
\midrule
Logistic ($p{=}1$) & MLP & $1.2\times10^{-4}$ & \textbf{$4.9\times10^{-5}$} & $7.4\times10^{-5}$ & $6.9\times10^{-5}$ & $4.7\times10^{-4}$ \\
 & B-spline KAN & $5.3\times10^{-5}$ & $4.2\times10^{-5}$ & $2.1\times10^{-5}$ & \textbf{$2.4\times10^{-5}$} & $5.1\times10^{-5}$ \\
 & Chebyshev KAN & $2.4\times10^{-5}$ & $2.9\times10^{-5}$ & $3.4\times10^{-5}$ & $3.2\times10^{-5}$ & \textbf{$2.4\times10^{-5}$} \\
\addlinespace
Singular IVP ($p{=}2$) & MLP & $4.8\times10^{-5}$ & $1.1\times10^{-4}$ & \textbf{$1.1\times10^{-5}$} & $1.8\times10^{-5}$ & $4.7\times10^{-5}$ \\
 & B-spline KAN & $1.9\times10^{-5}$ & $1.2\times10^{-5}$ & \textbf{$1.1\times10^{-5}$} & $1.8\times10^{-5}$ & $1.2\times10^{-5}$ \\
 & Chebyshev KAN & \textbf{$2.5\times10^{-5}$} & $4.7\times10^{-5}$ & $1.0\times10^{-5}$ & $3.4\times10^{-5}$ & $2.9\times10^{-5}$ \\
\addlinespace
Singular BVP ($p{=}2$) & MLP & \textbf{$3.7\times10^{-7}$} & $9.7\times10^{-6}$ & $1.4\times10^{-5}$ & $2.7\times10^{-5}$ & $1.5\times10^{-5}$ \\
 & B-spline KAN & $3.4\times10^{-7}$ & $1.5\times10^{-6}$ & $1.2\times10^{-6}$ & \textbf{$1.2\times10^{-7}$} & $1.5\times10^{-6}$ \\
 & Chebyshev KAN & $4.3\times10^{-7}$ & \textbf{$4.7\times10^{-7}$} & $5.0\times10^{-7}$ & $5.4\times10^{-7}$ & $1.3\times10^{-6}$ \\
\addlinespace
Falkner-Skan ($p{=}3$) & MLP & $5.2\times10^{-6}$ & $1.5\times10^{-6}$ & $4.9\times10^{-6}$ & \textbf{$1.6\times10^{-6}$} & $3.9\times10^{-6}$ \\
 & B-spline KAN & $6.8\times10^{-6}$ & $7.0\times10^{-6}$ & $4.1\times10^{-6}$ & $1.1\times10^{-5}$ & \textbf{$5.4\times10^{-7}$} \\
 & Chebyshev KAN & $2.9\times10^{-6}$ & $7.2\times10^{-6}$ & $1.5\times10^{-6}$ & \textbf{$4.8\times10^{-6}$} & $1.5\times10^{-6}$ \\
\addlinespace
Fourth-order IVP ($p{=}4$) & MLP & $8.3\times10^{-6}$ & $3.6\times10^{-6}$ & \textbf{$1.6\times10^{-6}$} & $1.8\times10^{-5}$ & $1.7\times10^{-4}$ \\
 & B-spline KAN & \textbf{$2.6\times10^{-6}$} & $5.5\times10^{-5}$ & $3.8\times10^{-7}$ & $3.6\times10^{-5}$ & $1.1\times10^{-2}$ \\
 & Chebyshev KAN & $1.1\times10^{-5}$ & $1.4\times10^{-5}$ & \textbf{$7.0\times10^{-7}$} & $7.3\times10^{-5}$ & $2.3\times10^{-5}$ \\
\bottomrule
\end{tabular}}
\ifdim\wd\fitbox>\linewidth\resizebox{\linewidth}{!}{\usebox\fitbox}\else\usebox\fitbox\fi
\end{table}
\Cref{tab:activations} and \cref{fig:activations} report the best-seed error of every activation. For the
MLP the activation is the node nonlinearity; for the KANs it is the residual base function added to the
learnable basis. No single activation is best for any architecture: the selected activation changes from
problem to problem (\cref{tab:main}). On the problems of order one and two, the MLP is the most sensitive to the choice. Its best and worst
activation differ by a factor of 9 to 73, because its entire nonlinearity is the activation. The KANs learn
their nonlinearity in the basis, and their best and worst activation differ by a factor of 1.4 to 13. At order three and four the B-spline KAN becomes sensitive as well: without the base
branch it does not fit the fourth-order problem at all, so the residual base function is essential there.

\subsection{Training dynamics}\label{sec:dynamics}
\begin{figure}[H]
\centering
\includegraphics[width=\linewidth]{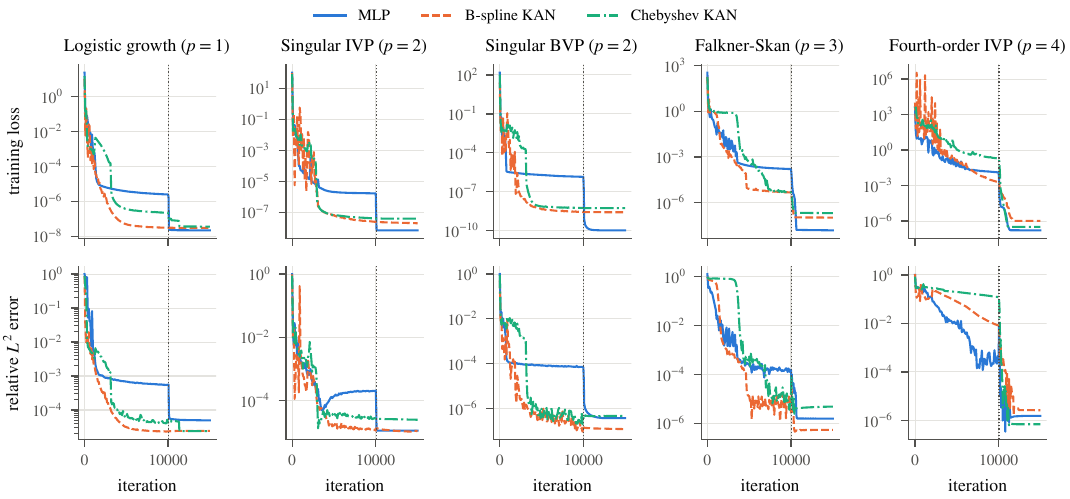}
\caption{Training loss (top) and relative $L^2$ test error (bottom) during training, best seed. The dotted
vertical line marks the switch from Adam to L-BFGS. The test error is recorded for plotting only and is not
used in training or model selection.}
\label{fig:training}
\end{figure}

\Cref{fig:training} shows the training loss and the test error during training; the vertical line marks the
switch from Adam to L-BFGS. On the problems of order one and two the KANs reach their final accuracy during
the Adam phase, and L-BFGS changes their loss by less than one order of magnitude; the MLP, in contrast,
stalls on a plateau during Adam and gains two to four orders of magnitude in loss, and roughly one to two in
error, from L-BFGS. On the third- and fourth-order problems, all three architectures depend on L-BFGS. Throughout,
the test error follows the training loss, which supports selecting runs by training loss.

\subsection{Parameter efficiency}\label{sec:scaling}
\begin{figure}[H]
\centering
\includegraphics[width=\linewidth]{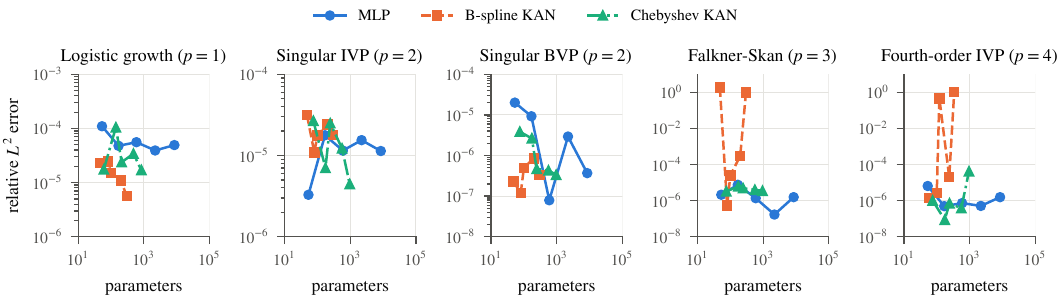}
\caption{Error versus number of trainable parameters (best of three seeds per size, logarithmic axes).}
\label{fig:scaling}
\end{figure}

\Cref{fig:scaling} varies the size of every architecture (MLP width 4 to 64; B-spline KAN width and grid;
Chebyshev KAN width and degree) under the same budget. The result is problem dependent and qualifies the
comparison of \cref{tab:main}, whose MLP is larger than it needs to be: on four of the five problems (all but
the singular BVP) MLPs with 53 to 169 parameters reach errors similar to, or lower than, the 8{,}513-parameter
baseline. On the logistic equation the B-spline KAN is more
accurate than the MLP at every size, and on the fourth-order problem the Chebyshev KAN attains the lowest error
of any network on that problem ($8.6\times10^{-8}$ with 175 parameters). On the singular problems P2 and P3 the three
families reach comparable errors at comparable sizes, and on the Falkner-Skan problem the MLP is the most
accurate. The B-spline KAN is unreliable on P4 and P5. Several of its configurations do not fit the solution, in
particular those that start from a finer grid, as \cref{rem:scale} predicts; they appear as isolated
high-error points in \cref{fig:scaling}.

\subsection{Grid extension}\label{sec:gridext}
\begin{figure}[H]
\centering
\includegraphics[width=\linewidth]{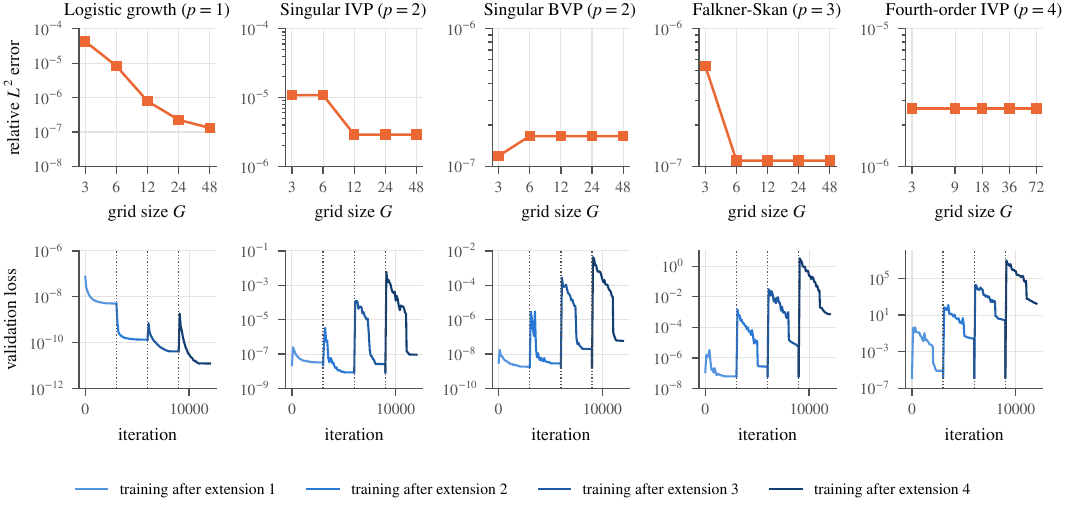}
\caption{Derivative-exact grid extension of the B-spline KAN, starting from the trained benchmark model
($G=3$; seed with the lowest final validation loss). Top: relative $L^2$ error after each stage versus grid size.
Bottom: validation loss (held-out collocation points) during the refinement stages; dotted lines mark the
extensions. The loss right after each extension equals the loss before it; the spikes that follow are caused by
the first optimizer steps on the finer grid.}
\label{fig:gridext}
\end{figure}

A B-spline KAN can be refined after training. With the usual least-squares transfer to a new, independently placed grid, which we tried first, the
physics loss jumped by 8
(Falkner-Skan) to 31 (fourth-order problem) orders of magnitude after the extensions for $p\ge2$, although the
solution values barely changed, and the refined networks did not improve. With the knot-insertion transfer of \cref{rem:nested} the refined network is the old
network, and the loss is unchanged to four significant digits. \Cref{fig:gridext} and \cref{tab:gridext} start
from the trained benchmark models ($G=3$) and double the grid four times, training after each extension with
a small learning rate and selecting checkpoints on held-out collocation points.
In all 60 extensions (five problems, three seeds, four stages) the validation loss immediately after the
extension equals the loss before it: the transfer is exact.

The benefit of the added resolution depends on the differential order. On the logistic equation the error
decreases at every stage, by a factor of \gridGainLogistic{} in total for the selected seed (160 to 690 for the
three seeds), to $1.3\times10^{-7}$ at $G=\gridFinalLogistic$. That is more than two orders of magnitude below
the best MLP. On the singular IVP and the Falkner-Skan problem the error decreases by factors of
\gridGainSingularIVP{} and \gridGainFalknerSkan{} and then levels off, reaching $1.1\times10^{-7}$ on the
Falkner-Skan problem, fourteen times below the best MLP. On the singular BVP the coarse model is already at the
$10^{-7}$ level and the extension changes the error by a factor of \gridGainSingularBVP{}, and on the fourth-order
problem it brings no improvement (factor \gridGainFourthOrder{}); for two of the three seeds of P5 the selected
checkpoint even has a larger error, although a smaller residual, than the unrefined network. The bottom row of
\cref{fig:gridext} shows the reason: the extension itself leaves the loss unchanged, but the first Adam steps on the
finer grid raise it sharply, and the more so the higher the differential order, because a small change of the
coefficients of a fine spline changes its $p$-th derivative strongly (\cref{rem:scale}). Refinement helps most at
low order, and a gentler optimizer in the refinement stages (for example L-BFGS alone) is a
natural next step for $p\ge3$. The refinement stages add training iterations, so these gains are not part of the
equal-budget comparison of \cref{tab:main}. Neither the MLP nor the Chebyshev KAN has an equivalent operation:
enlarging them means training a new network.

\begin{table}[t]
\centering
\caption{Grid extension of the B-spline KAN, starting from the trained benchmark model (stage 0, $G=3$). Each stage doubles the grid by exact knot insertion and continues training; relative $L^2$ error of the seed with the lowest final validation loss, grid size in parentheses.}
\label{tab:gridext}
\small
\setlength{\tabcolsep}{3pt}
\ifdefined\fitbox\else\newsavebox\fitbox\fi
\sbox\fitbox{\begin{tabular}{lccccc}
\toprule
Problem & Stage 0 & Stage 1 & Stage 2 & Stage 3 & Stage 4 \\
\midrule
Logistic ($p{=}1$) & $4.3\times10^{-5}$\,{\scriptsize(3)} & $8.3\times10^{-6}$\,{\scriptsize(6)} & $7.9\times10^{-7}$\,{\scriptsize(12)} & $2.3\times10^{-7}$\,{\scriptsize(24)} & $1.3\times10^{-7}$\,{\scriptsize(48)} \\
Singular IVP ($p{=}2$) & $1.1\times10^{-5}$\,{\scriptsize(3)} & $1.1\times10^{-5}$\,{\scriptsize(6)} & $2.9\times10^{-6}$\,{\scriptsize(12)} & $2.9\times10^{-6}$\,{\scriptsize(24)} & $2.9\times10^{-6}$\,{\scriptsize(48)} \\
Singular BVP ($p{=}2$) & $1.2\times10^{-7}$\,{\scriptsize(3)} & $1.7\times10^{-7}$\,{\scriptsize(6)} & $1.7\times10^{-7}$\,{\scriptsize(12)} & $1.7\times10^{-7}$\,{\scriptsize(24)} & $1.7\times10^{-7}$\,{\scriptsize(48)} \\
Falkner-Skan ($p{=}3$) & $5.4\times10^{-7}$\,{\scriptsize(3)} & $1.1\times10^{-7}$\,{\scriptsize(6)} & $1.1\times10^{-7}$\,{\scriptsize(12)} & $1.1\times10^{-7}$\,{\scriptsize(24)} & $1.1\times10^{-7}$\,{\scriptsize(48)} \\
Fourth-order IVP ($p{=}4$) & $2.6\times10^{-6}$\,{\scriptsize(3)} & $2.6\times10^{-6}$\,{\scriptsize(9)} & $2.6\times10^{-6}$\,{\scriptsize(18)} & $2.6\times10^{-6}$\,{\scriptsize(36)} & $2.6\times10^{-6}$\,{\scriptsize(72)} \\
\bottomrule
\end{tabular}}
\ifdim\wd\fitbox>\linewidth\resizebox{\linewidth}{!}{\usebox\fitbox}\else\usebox\fitbox\fi
\end{table}

\subsection{Interpretability}\label{sec:interp}
\begin{figure}[H]
\centering
\includegraphics[width=\linewidth]{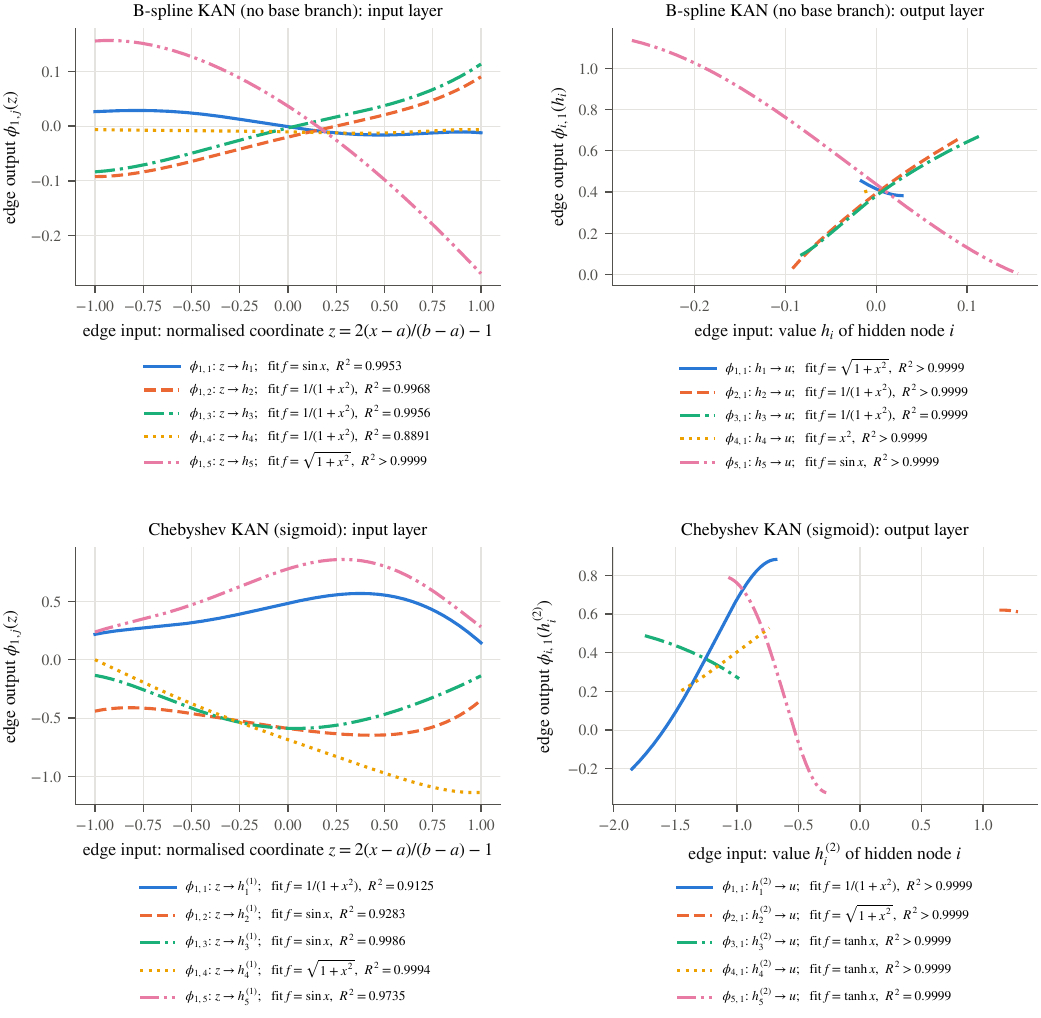}
\caption{Learned edge functions of the two KANs trained on the Falkner-Skan problem (top: B-spline KAN, bottom:
Chebyshev KAN; the selected base function is given in each title). Left: input layer, one curve per edge
$\phi_{1,j}$ from the normalised coordinate $z$ to hidden node $j$. Right: output layer, one curve per edge
$\phi_{i,1}$ from hidden node $i$ to the output $u$, plotted over the values $h_i$ that node $i$ takes on the domain.
Each edge has its own color and line style. The legend names the edge, the function $f$ for which
$a\,f(bx+c)+d$ fits the curve best, and the coefficient of determination $R^2$ of that fit.}
\label{fig:edges}
\end{figure}

In a KAN every edge carries a function of one variable, so the trained network can be inspected by plotting
these functions. We write $\phi_{i,j}$ for the function on the edge from node $i$ of one layer to node $j$ of the
next. \Cref{fig:edges} shows the trained KANs of the Falkner-Skan problem. In the input layer there is one input
node, the normalized coordinate $z$, so the five edges are $\phi_{1,1},\dots,\phi_{1,5}$, and $\phi_{1,j}$ sends $z$
to hidden node $j$. In the output layer the five edges $\phi_{1,1},\dots,\phi_{5,1}$ send the value $h_i$ of hidden
node $i$ to the output $u$. The Chebyshev KAN has two hidden layers ($h^{(1)}$ and $h^{(2)}$); its 25 edges between them
are not shown. Each curve is drawn over the range of values that its edge actually receives on the domain, which is
why the output-layer curves have different lengths.

To see whether an edge is close to a simple formula, we fitted each curve by $a\,f(bx+c)+d$, where $f$ is one of
$x$, $x^2$, $x^3$, $\exp x$, $\sin x$, $\tanh x$, $1/(1+x^2)$ and $\sqrt{1+x^2}$, and $a,b,c,d$ are fitted
constants. The legend of \cref{fig:edges} gives, for every edge, the function $f$ with the best fit and its
coefficient of determination $R^2$ ($R^2=1$ is a perfect fit); all fitted constants are listed in the files that
accompany the code. Over the five problems, between 65\% and 95\% of the edges are matched with $R^2>0.99$, and the
median $R^2$ is at least $0.998$ on every problem. Edges with a lower $R^2$ are nearly constant (for example
$\phi_{1,4}$ of the B-spline KAN), so their small variation is harder to match but contributes little. Such fits are
the first step of the symbolic regression procedure of \citep{liu2024kan}. An MLP offers no equivalent view: none of
its 8{,}513 weights corresponds to a function of the input.

\subsection{Robustness and cost}\label{sec:robust}
\begin{table}[t]
\centering
\caption{Behaviour across all runs (every activation and seed). \emph{Finished}: runs whose final model is finite (the training loop rolls back any non-finite step, see Sec.~\ref{sec:training}). \emph{Not fitted}: finished runs whose relative $L^2$ error exceeds $10^{-2}$. Median and worst relative $L^2$ error over all runs.}
\label{tab:robust}
\small
\setlength{\tabcolsep}{3pt}
\ifdefined\fitbox\else\newsavebox\fitbox\fi
\sbox\fitbox{\begin{tabular}{llcccc}
\toprule
Problem & Architecture & Finished & Not fitted & Median rel.\ $L^2$ & Worst rel.\ $L^2$ \\
\midrule
Logistic ($p{=}1$) & MLP & 15/15 & 0 & $3.3\times10^{-4}$ & $5.4\times10^{-4}$ \\
 & B-spline KAN & 15/15 & 0 & $5.8\times10^{-5}$ & $1.5\times10^{-4}$ \\
 & Chebyshev KAN & 15/15 & 0 & $3.2\times10^{-5}$ & $4.9\times10^{-5}$ \\
\addlinespace
Singular IVP ($p{=}2$) & MLP & 15/15 & 0 & $7.8\times10^{-5}$ & $7.2\times10^{-4}$ \\
 & B-spline KAN & 15/15 & 0 & $1.8\times10^{-5}$ & $7.5\times10^{-5}$ \\
 & Chebyshev KAN & 15/15 & 0 & $3.6\times10^{-5}$ & $8.5\times10^{-5}$ \\
\addlinespace
Singular BVP ($p{=}2$) & MLP & 15/15 & 0 & $1.6\times10^{-5}$ & $7.1\times10^{-5}$ \\
 & B-spline KAN & 15/15 & 0 & $1.4\times10^{-6}$ & $5.1\times10^{-6}$ \\
 & Chebyshev KAN & 15/15 & 0 & $6.6\times10^{-7}$ & $2.8\times10^{-6}$ \\
\addlinespace
Falkner-Skan ($p{=}3$) & MLP & 15/15 & 0 & $5.1\times10^{-6}$ & $2.9\times10^{-5}$ \\
 & B-spline KAN & 15/15 & 0 & $5.8\times10^{-6}$ & $2.4\times10^{-5}$ \\
 & Chebyshev KAN & 15/15 & 3 & $4.8\times10^{-6}$ & $1.8\times10^{0}$ \\
\addlinespace
Fourth-order IVP ($p{=}4$) & MLP & 15/15 & 0 & $9.0\times10^{-6}$ & $5.2\times10^{-4}$ \\
 & B-spline KAN & 15/15 & 5 & $2.6\times10^{-4}$ & $2.7\times10^{-1}$ \\
 & Chebyshev KAN & 15/15 & 0 & $1.9\times10^{-5}$ & $2.8\times10^{-4}$ \\
\bottomrule
\end{tabular}}
\ifdim\wd\fitbox>\linewidth\resizebox{\linewidth}{!}{\usebox\fitbox}\else\usebox\fitbox\fi
\end{table}
\Cref{tab:robust} reports every run. All \nRuns{} runs finished with a finite model, and the rollback of
\cref{alg:training} was never triggered. \nUnfit{} runs nevertheless did not fit the solution (relative $L^2$
error above $10^{-2}$): three Chebyshev KAN runs on the Falkner-Skan problem, all with the same seed, and
five B-spline KAN runs on the fourth-order problem, including all three runs without a base branch. Every MLP run fitted
its problem. So the KANs are more accurate at their best but less reliable at order three and four. For this
reason, the median and the worst run are reported next to the best one. Per run, the B-spline KAN is as fast as
the MLP on the problems of order one and two but 1.3 to 1.9 times slower on P4 and P5, because the Cox-de Boor
recursion is differentiated $p$ times; the Chebyshev KAN is 1.3 to 1.7 times faster than the MLP throughout.

\section{Discussion}\label{sec:discussion}
\paragraph{Parameter count.} The KANs in this study have \params{FalknerSkan}{BSpline} and
\params{FalknerSkan}{Cheby} parameters, against \params{FalknerSkan}{MLP} for the MLP.
Their best runs are as accurate as, or up to three times more accurate than, the baseline MLP on every
problem. The size sweep shows, however, that the baseline MLP is itself oversized: much smaller MLPs are
competitive, so the parameter advantage of KANs over a well-sized MLP is real on the logistic and fourth-order
problems but not general. A parameter count is also not a runtime: at order three and four the B-spline KAN is the
slowest per run.

\paragraph{High-order residuals favour coarse bases.} For physics-informed KANs, \cref{rem:scale} has a
practical consequence: a basis that is fine enough to represent the solution accurately is too stiff to be
optimised from a random start when the residual contains third or fourth derivatives. The Chebyshev KAN behaves
analogously with the polynomial degree. The basis has to start coarse, and a later refinement has to preserve derivatives as well as values
(\cref{rem:nested}).

\paragraph{Limitations.} All problems are scalar ODEs on bounded intervals; loss weights are fixed rather
than adaptive \citep{wang2021gradient}; the KAN sizes were chosen by a short pilot on two problems; and the main
comparison reports the best of fifteen runs per architecture, which favors the architectures with the larger
spread between runs. The median and worst runs in \cref{tab:robust} should be read alongside it. Extending the comparison to systems of ODEs (such as
the coupled nanofluid model of \citep{venkatachalapathy2023}) and to partial differential equations is
future work.

\section{Conclusion}\label{sec:conclusion}
In this study, a comprehensive benchmarking of multilayer perceptrons, B-spline Kolmogorov-Arnold Networks, and Chebyshev KANs was conducted for physics-informed solutions of nonlinear ordinary differential equations spanning orders one through four. Operating under an identical, standardized training budget and protocol, all \nRuns{} training runs achieved numerical stability without divergence. The findings demonstrate that KAN architectures match or exceed the accuracy of baseline MLPs while utilizing one to two orders of magnitude fewer parameters. However, parameter scaling evaluations reveal that compact MLPs remain highly competitive, and standard MLPs exhibit greater robustness, successfully fitting every realization across higher-order problems where KAN models occasionally faltered.

Beyond raw parameter efficiency, KANs provide unique structural benefits for scientific computing. Their edge-based design facilitates rapid convergence during the initial Adam phase on low-order problems, and the explicit univariate functions mapped onto every edge allow for straightforward inspection and closed-form analytical matching. Moreover, replacing traditional grid extension with derivative-preserving knot insertion effectively bypasses residual distortion, yielding substantial error reductions at lower orders (over two orders of magnitude on the logistic equation). Ultimately, while KANs offer exceptional accuracy and interpretability, high-order differential equations inherently require careful initialization with coarse bases to ensure stable optimization.

\appendix
\section{Pointwise errors}\label{app:pointwise}
\begin{table}[p]
\centering
\caption{Pointwise absolute error $|y_{\mathrm{ref}}-\hat y|$ of the best run of each architecture at equally spaced points of each domain, together with the exact or reference solution.}
\label{tab:pointwise}
\small
\setlength{\tabcolsep}{5pt}
\ifdefined\fitbox\else\newsavebox\fitbox\fi
\sbox\fitbox{\begin{tabular}{lccccc}
\toprule
Problem & $x$ & $y_{\mathrm{ref}}(x)$ & MLP & B-spline KAN & Chebyshev KAN \\
\midrule
Logistic ($p{=}1$) & 0.00 & 0.500000 & $1.7\times10^{-9}$ & $2.4\times10^{-9}$ & $1.5\times10^{-8}$ \\
 & 1.00 & 0.731059 & $4.1\times10^{-5}$ & $7.8\times10^{-6}$ & $2.2\times10^{-5}$ \\
 & 2.00 & 0.880797 & $6.5\times10^{-5}$ & $3.0\times10^{-5}$ & $1.0\times10^{-5}$ \\
 & 3.00 & 0.952574 & $6.1\times10^{-5}$ & $1.9\times10^{-6}$ & $1.6\times10^{-5}$ \\
 & 4.00 & 0.982014 & $1.2\times10^{-5}$ & $7.9\times10^{-6}$ & $2.0\times10^{-5}$ \\
 & 5.00 & 0.993307 & $2.3\times10^{-5}$ & $3.8\times10^{-6}$ & $9.3\times10^{-6}$ \\
\addlinespace
Singular IVP ($p{=}2$) & 0.00 & 1.000000 & $3.1\times10^{-9}$ & $3.7\times10^{-8}$ & $5.0\times10^{-9}$ \\
 & 0.20 & 0.980581 & $1.5\times10^{-5}$ & $1.3\times10^{-5}$ & $3.3\times10^{-5}$ \\
 & 0.40 & 0.928477 & $1.3\times10^{-5}$ & $1.3\times10^{-5}$ & $2.8\times10^{-5}$ \\
 & 0.60 & 0.857493 & $9.8\times10^{-6}$ & $8.4\times10^{-6}$ & $1.7\times10^{-5}$ \\
 & 0.80 & 0.780869 & $4.9\times10^{-6}$ & $3.3\times10^{-6}$ & $6.1\times10^{-6}$ \\
 & 1.00 & 0.707107 & $7.4\times10^{-8}$ & $1.0\times10^{-6}$ & $3.3\times10^{-6}$ \\
\addlinespace
Singular BVP ($p{=}2$) & 0.00 & 0.954135 & $8.6\times10^{-7}$ & $6.0\times10^{-7}$ & $3.8\times10^{-7}$ \\
 & 0.20 & 0.955950 & $6.9\times10^{-8}$ & $9.1\times10^{-8}$ & $5.2\times10^{-7}$ \\
 & 0.40 & 0.961403 & $5.2\times10^{-7}$ & $1.4\times10^{-7}$ & $3.1\times10^{-7}$ \\
 & 0.60 & 0.970526 & $1.2\times10^{-7}$ & $2.0\times10^{-8}$ & $2.8\times10^{-7}$ \\
 & 0.80 & 0.983369 & $1.2\times10^{-7}$ & $1.3\times10^{-7}$ & $7.4\times10^{-8}$ \\
 & 1.00 & 1.000000 & $1.5\times10^{-9}$ & $2.6\times10^{-9}$ & $3.3\times10^{-9}$ \\
\addlinespace
Falkner-Skan ($p{=}3$) & 0.00 & 0.916674 & $1.6\times10^{-8}$ & $1.3\times10^{-10}$ & $2.9\times10^{-8}$ \\
 & 0.60 & 1.113635 & $4.8\times10^{-7}$ & $5.0\times10^{-8}$ & $1.2\times10^{-6}$ \\
 & 1.20 & 1.562807 & $1.6\times10^{-6}$ & $7.6\times10^{-7}$ & $5.8\times10^{-6}$ \\
 & 1.80 & 2.118838 & $3.1\times10^{-6}$ & $1.6\times10^{-6}$ & $1.1\times10^{-5}$ \\
 & 2.40 & 2.709473 & $4.9\times10^{-6}$ & $1.5\times10^{-6}$ & $1.5\times10^{-5}$ \\
 & 3.00 & 3.308343 & $5.4\times10^{-6}$ & $1.1\times10^{-6}$ & $1.5\times10^{-5}$ \\
\addlinespace
Fourth-order IVP ($p{=}4$) & 0.00 & 1.000000 & $4.7\times10^{-7}$ & $7.7\times10^{-7}$ & $1.8\times10^{-7}$ \\
 & 0.20 & 0.833333 & $4.9\times10^{-7}$ & $8.2\times10^{-7}$ & $1.8\times10^{-7}$ \\
 & 0.40 & 0.714286 & $5.9\times10^{-7}$ & $9.9\times10^{-7}$ & $2.1\times10^{-7}$ \\
 & 0.60 & 0.625000 & $8.3\times10^{-7}$ & $1.5\times10^{-6}$ & $3.4\times10^{-7}$ \\
 & 0.80 & 0.555556 & $1.4\times10^{-6}$ & $2.5\times10^{-6}$ & $6.4\times10^{-7}$ \\
 & 1.00 & 0.500000 & $2.5\times10^{-6}$ & $4.2\times10^{-6}$ & $1.2\times10^{-6}$ \\
\bottomrule
\end{tabular}}
\ifdim\wd\fitbox>\linewidth\resizebox{\linewidth}{!}{\usebox\fitbox}\else\usebox\fitbox\fi
\end{table}
\clearpage

\section*{Declarations}
\noindent\textbf{Funding} Not applicable.\\
\textbf{Conflict of interest} The author declares no conflict of interest.\\
\textbf{Data and code availability} All code, raw results of every seed, and the scripts that generate every
table and figure of this paper are provided in the accompanying repository.

\bibliographystyle{plain}
\bibliography{p1_ref}

\end{document}